\documentclass[a4paper,12pt,reqno]{article}

\usepackage{amsmath}
\usepackage{amsthm}
\usepackage{amssymb}
\usepackage{latexsym}
\usepackage{graphicx}
\usepackage{stmaryrd}
\usepackage{mathrsfs}
\usepackage{enumerate}
\usepackage{cases}

\newcommand{\glalign}[2]{\lower.6ex\vbox{
\baselineskip\lineskip\ialign{$#1\hfil##\hfil$\crcr#2\crcr=\crcr}}}

\newcommand{\R}{\mathbb{R}}

\newcommand{\pt}{\partial_t}
\newcommand{\del}{\partial}

\newcommand{\delt}{\partial_{t}}

\renewcommand{\div}{\mbox{\rm div}\,}

\newcommand{\trans}{{}^\top}

\def\eqn#1$$#2$${\begin{equation}\label#1#2\end{equation}}

\numberwithin{equation}{section}

\newtheorem{defi}{Definition}[section]
\newtheorem{thm}[defi]{Theorem}
\newtheorem{cor}[defi]{Corollary}
\newtheorem{prop}[defi]{Proposition}
\newtheorem{lem}[defi]{Lemma}
\newtheorem{rem}[defi]{Remark}

\def\eqn#1$$#2$${\begin{equation}\label#1#2\end{equation}}

\numberwithin{equation}{section}
\numberwithin{equation}{section}
\allowdisplaybreaks[4]
\begin{document}

\title{
\bf \large 
Asymptotic profile for diffusion wave terms of the compressible Navier-Stokes-Korteweg system}
\author{
Takayuki KOBAYASHI\footnote{Osaka University, 1-3, Machikaneyamacho, Toyonakashi, 560-8531, JAPAN}\\
Masashi MISAWA\footnote{Kumamoto University, 2-39-1, Kurokami,  Chuo-ku,  Kumamoto, 
860-8555, JAPAN}\\
Kazuyuki TSUDA\footnote{Osaka University, 1-3, Machikaneyamacho, Toyonakashi, 560-8531, JAPAN, corresponding author,  \qquad mail: jtsuda@sigmath.es.osaka-u.ac.jp}\\
}
\date{}
\maketitle
\begin{abstract} 
 Asymptotic profile for diffusion wave terms of solutions to the compressible Navier-Stokes-Korteweg system is studied on $\mathbb{R}^2$.  
The diffusion wave with time decay estimate is studied by Hoff and Zumbrun (1995, 1997), Kobayashi and Shibata (2002) and Kobayashi and Tsuda (2018) for the compressible Navier-Stokes system and the compressible Navier-Stokes-Korteweg system.
Our main assertion in this paper is that, for some initial conditions given by the Hardy space,
asymptotic behaviors in space-time $L^2$ of the diffusion wave parts
are essentially different between density and the potential flow part of the momentum.
Even though measuring by $L^2$ on space,
a decay of the potential flow part is slower than that of the Stokes flow part of the momentum.  
The proof is based on a modified version of Morawetz's energy estimate,
and the Fefferman-Stein inequality on the duality between the Hardy space
and functions of bounded mean oscillation.
\end{abstract}



%

\section{Introduction}

We study asymptotic behavior of solutions to the following compressible Navier-Stokes-Korteweg system in $\mathbb{R}^2$,
called ``CNSK'':
\begin{numcases}
{}
\partial_{t}\rho +\div M=0,\nonumber\\
\partial_{t}M +\div \Big(\frac{M \otimes M}{\rho}\Big)+\nabla P(\rho)=\div \Big({\mathcal S}(\frac{M}{\rho})+{\mathcal K}(\rho)\Big),\label{CNSK}\\
\rho(x,0) =\rho_0, \ \ M(x,0)=M_0. \nonumber
\end{numcases}
\noindent Here $\rho=\rho(x,t)$ and $M=(M_{1}(x,t),M_{2}(x,t))$ are the unknown density and momentum, respectively, at time $t\in \mathbb{R}^{+}$ and 
position $x\in\mathbb{R}^2$; $\rho_0=\rho_0(x)$ and $M_0=M_0(x)$ are given initial data; 
${\mathcal S}$ and ${\mathcal K}$ denote the viscous stress tensor and the Korteweg stress tensor, respectively,
given by 
\begin{eqnarray}
\left\{
\begin{array}{ll}
({\mathcal S})_{i,j}(\frac{M}{\rho}) =\Big(\mu' \div \frac{M}{\rho} \Big)\delta_{i,j} +2\mu d_{ij}\Big(\frac{M}{\rho}\Big), \\
({\mathcal K})_{i,j}(\rho) =\frac{\kappa}{2}(\Delta \rho^2 -|\nabla \rho|^2)\delta_{i,j}- \kappa \frac{\del \rho}{\del x_i}\frac{\del \rho}{\del x_j},
\end{array}
\right.\label{Korteweg tensor}
\end{eqnarray}
where $d_{ij}\Big(\frac{M}{\rho}\Big)=\frac{1}{2}\left(\frac{\del }{\del x_i}\Big(\frac{M}{\rho}\Big)_j+\frac{\del}{\del x_j}\Big(\frac{M}{\rho}\Big)_i\right)$; $\mu$ and $\mu'$ are the viscosity coefficients,
supposed to be constants satisfying 
$$
\mu>0, \ \ \ 
\mu+\mu'\geq 0. 
$$
$\kappa$ is
the capillary constant satisfying
$\kappa \geq 0$.
If $\kappa=0$ in the Korteweg tensor, the usual compressible Navier-Stokes equation,
abbreviated to ``CNS'', appears;
$P=P(\rho)$ is the pressure,
assumed to be a smooth function of $\rho$ satisfying  
$P'(\rho_*)>0$,
where $\rho_*$ is a given positive constant and $(\rho_*,0)$ is
a given constant state for \eqref{CNSK}. 
We consider solutions to \eqref{CNSK} around the constant state. 

\eqref{CNSK} is the system of equations of motion of liquid-vapor type two phase flow
with phase transition in a compressible fluid, similarly as in \cite{de-bu}.
To describe the phase transition, this model use the diffusive interface.
Hence the phase boundary is regarded as a narrow transition layer and
change of the density
prescribes fluid state.    
Due to the diffusive interface, it is enough to consider one set of equations and  a single spatial domain and difficulty of topological change of interface does not occur. 

For derivation of \eqref{CNSK}, 
Van der Waals \cite{Van der Waals} suggests that a phase transition boundary is regarded as a thin transition zone, i.e, diffusive interface caused by a steep gradient of the density.  Based on his idea, Korteweg \cite{Korteweg}  modifies the stress tensor  of the Navier-Stokes equation to that including the term $\nabla\rho \otimes \nabla \rho$. 
Dunn and Serrin \cite{Dunn and Serrin} generalize the Korteweg's work
and strictly provide the system $(\ref{CNSK})$ with (\ref{Korteweg tensor}).
In their recent works, Heida and M\'{a}lek \cite{Heida and malek} 
derive \eqref{CNSK} by  the entropy production method.  

We will focus on the ``diffusion wave'' which stems from hyperbolic and parabolic aspects
of the system. The diffusion wave is given by convolution between the heat kernel and fundamental solution to the wave equation. The importance of diffusion wave for problems in one dimensional case was first recognized by Liu \cite{Liu} for the study of stability of shock waves for viscous conservation laws.  
The multi-dimensional diffusion wave with time decay estimate of solutions is studied for CNS by Hoff and Zumbrun \cite{Hoff-Zumbrun1, Hoff-Zumbrun2} and Kobayashi and Shibata \cite{Kobayashi-Shibata}, and for the viscoelastic equation
on $\mathbb{R}^n$ $(n \geq 2)$ by Shibata \cite{Shibata-v}.
Let $u=\trans(\rho-\rho_*, M)$ be a solution to CNS and set  $E:=\|u_0\|_{H^{s+\ell}\cap L^1}$, where $u_0= \trans(\rho_0-\rho_*, M_0)= \trans(\rho(0)-\rho_*, M(0))$, $s$ is an integer part of $n/2$ and $\ell$ is integer satisfying $\ell \geq 3$.
Then,  it is shown in \cite{Hoff-Zumbrun1, Hoff-Zumbrun2, Kobayashi-Shibata}
that the linear parts decays faster than nonlinear parts in the Duhamel formula and the asymptotic behavior in $L^p (\mathbb{R}^n)$ $(p >2, n \geq 2)$ of solutions is presented as 
\begin{eqnarray}
u(t) \sim 
\underbrace{
\begin{pmatrix}
0 \\
{\cal K}_{\nu}(t)\ast M_{0,in}
\end{pmatrix}
+ \begin{pmatrix}
\rho(t)- \rho_*\\
M(t)-{\cal K}_{\nu}(t)\ast M_{0,in}
\end{pmatrix}
}_{\mbox{solutions to linearized system}}
+\underbrace{
 \cdots}_{\mbox{nonlinear parts}}
\ \ \mbox{in}  \ \ L^p(\mathbb{R}^n) \label{behavior-cns}
\end{eqnarray}
as $t$ goes to infinity.
Here, the notation $u (t) \sim f (t)$ in $L^p (\mathbb{R}^n)$
is defined as
$$
\limsup_{t \rightarrow \infty}
\|u (t) - f (t)\|_{L^p (\mathbb{R}^n)}
\le C
$$
for a positive number $C$ independent of $t$,
the similar notation will be used hereafter.
${\cal K}_{\nu}={\cal K}_{\nu}(t,x)$ is
the standard heat kernel and $M_{0,in}$ is
a divergence-free part of $M_0$,
given by
$$
{\cal K}_{\nu}=\mathcal{F}^{-1}(e^{-\nu|\xi|^2 t}), \ \ M_{0,in}=\mathcal{F}^{-1}\Big\{\Big(I_n-\frac{\xi\trans{\xi}}{|\xi|^2}\Big)M_0\Big\}.
$$
More precisely, it holds that 
\begin{equation*}
\|u(t)\|_{L^p (\mathbb{R}^n)}\leq 
C(1+t)^{-\frac{n}{2}(1-\frac{1}{p})},  \ \ 
\|u(t)-G(t)\ast u_0 \|_{L^p (\mathbb{R}^n)}\leq 
C(1+t)^{-\frac{n}{2}(1-\frac{1}{p})-\frac{1}{2}}
\end{equation*}
and 
\begin{eqnarray*}
\Big\|
\begin{pmatrix}
\rho(t) - \rho_*
\\
M(t) - {\cal K}_\nu (t) \ast M_{0,in}
\end{pmatrix}
\Big\|_{L^p (\mathbb{R}^n)}  \leq C(1+t)^{-\frac{n}{2}(1-\frac{1}{p})+(\frac{n-1}{4})(\frac{2}{p}-1)} 
\end{eqnarray*}
for $t>0$, where $G=G(t)$ is the Green function of linearized CNS and $(\frac{n-1}{4})(\frac{2}{p}-1)\leq 0$ when $2\leq p \leq \infty$.   
Note that ${\cal K}_{\nu} (t) \ast M_{0,in}$ and $M(t)-{\cal K}_{\nu} (t) \ast M_{0,in}$
are the Stokes flow and  potential flow parts of $M$, respectively,
in the Helmholtz decomposition.
$\rho(t)- \rho_*$ and $M(t)-{\cal K}_{\nu} (t) \ast M_{0,in}$ are given
by the Green matrix of the linearized system,
which consists of the convolution with the Green functions 
of the diffusion and the wave equations
and thus, 
are called as the diffusion wave part. In addition, when $p=2$ the behaviors of both of $\rho(t)- \rho_*$ and $M(t)-{\cal K}_{\nu} (t) \ast M_{0,in}$ coincide with 
${\cal K}_{\nu}(t)\ast M_{0,in}$ as the parabolic type decay rate.
Kobayashi and Tsuda study the diffusion wave property for \eqref{CNSK} in \cite{Kobayashi-Tsuda}.

In this paper we consider the linearized system for \eqref{CNSK}.
Under some initial conditions given by the Hardy space ${\cal H}^1$(defined below),
we show some space-time $L^2$estimates
for the density and the Stokes flow part of the momentum.
The potential flow part of the momentum is also shown to grow
at the rate of logarithmic order in spatial-time $L^2$ norm.
See the precise initial condition given by the Hardy space below.
%
Here we assume a stronger initial condition  by ${\cal H}^1$ for density
than that by $L^1$,
in contrast to \cite{Kobayashi-Tsuda},
and thus, our results may show a gain of regularity by the Hardy space in the decay estimates.
Such a gain is also obtained for the heat equations.  In fact, we consider the following Cauchy problem:  
\begin{eqnarray}
\left\{
\begin{array}{ll}
\del_t u- \Delta u =0, \label{heat-eq}  \ \ \mbox{in}  \ \ (0, \infty) \times \mathbb{R}^2, \label{heat-eq}\\
u=u_0 \ \ \mbox{on}  \ \ \{t=0\} \times \mathbb{R}^2. 
\end{array}
\right.
\end{eqnarray}
It is well known that the solution $u$ to \eqref{heat-eq} satisfies the estimate
$$
\displaystyle\int_0^\infty \|u(t)\|_{L^2(\mathbb{R}^2)}^2 dt  \leq C\|u_0\|_{{\cal H}^1}^2
$$
for $u_0 \in {\cal H}^1$, while ${\cal H}^1 \subset L^1$ and for $u_0 \in L^1$ the estimate 
$$
\displaystyle\int_0^\infty \|u(t)\|_{L^2(\mathbb{R}^2)}^2 dt  \leq C\|u_0\|_{L^1}^2 
$$
generally does not hold.  
This fact shows a suitable gain of regularity by the
Hardy space and motivates us to introduce the Hardy space in this present paper.  
The nonlinear terms are expected
to decay in faster than the linear terms in the Duhamel formula,
similarly as in \cite{Hoff-Zumbrun1, Kobayashi-Tsuda}.
As a consequence,
the leading terms of the asymptotic expansion of the solution $u$ for \eqref{CNSK} is given by 
\begin{eqnarray}
u(t) \sim \underbrace{
\begin{pmatrix}
0 \\
M(t)-{\cal K}_{\nu}(t)\ast M_{0,in}
\end{pmatrix}
+ \begin{pmatrix}
\rho(t)- \rho_*\\
{\cal K}_{\nu}(t)\ast M_{0,in}
\end{pmatrix}
}_{\mbox{solutions to linearized system}}
+\underbrace{
 \cdots}_{\mbox{nonlinear parts}}
\mbox{in}\, L^2(0,\infty; \mathbb{R}^2). \label{behavior}
\end{eqnarray}
Precisely, the following estimates hold true for the solutions to the linearized CNSK :
$$
\displaystyle
\int_0^t \|M(\tau)-{\cal K}_{\nu}(\tau)\ast M_{0,in}\|_{L^2 (\mathbb{R}^2)}^2 d \tau
\sim \log t  \ \ \mbox{as} \ \ t \rightarrow \infty,  
$$
$$
\|\rho- \rho_*\|_{L^2 (0,\infty; \, L^2 (\mathbb{R}^2))}
+
\|{\cal K}_{\nu}\ast M_{0,in}\|_{L^2 (0,\infty; \, L^2 (\mathbb{R}^2))} < \infty.  
$$
The behaviors above of the diffusion wave parts $\rho(t)- \rho_*$
and $M(t)-{\cal K}_{\nu} (t) \ast M_{0,in}$  are clearly different from \eqref{behavior-cns}.
Even though measuring by $L^2$ on space,
$M(t)-{\cal K}_{\nu} (t) \ast M_{0,in}$ decays slower than the Stokes flow part of $M$.  
By the dependence on $\kappa$ of constants,
the above estimate \eqref{behavior} also holds true for CNS
(Theorems \ref{time-spatial est} and \ref{time-spatial est-CNS}).
We also obtain a decay rate of $L^2$ norm of the density
(Theorem \ref{time-decay est-CNSK}).
Furthermore, if $M_0 \in {\cal H}^1$,
the space-time $L^2$ boundedness is obtained for
$M(t)-{\cal K}_{\nu} (t) \ast M_{0,in}$, $\rho(t)- \rho_*$ and ${\cal K}_{\nu}(t)\ast M_{0,in}$.

The proofs of main results are based on the Morawetz type energy estimates for the linearized system.
We show that the diffusion wave part of the density $\rho(t)- \rho_*$
is bounded in
space-time $L^2$.
We rewrite \eqref{CNSK} to some linear doubly dispersion equation for $\rho$
and apply a modified version of Morawetz's energy estimate.
A preliminary function is introduced in the Morawetz estimate (see \eqref{w-okikata} below),
which is defined by use of a doubly Laplace type equation.
The existence of solution to the linear doubly Laplace type equation
is shown by use of the linear theory on ${\cal H}^1$,
which may be of its own interest.
Through the preliminary function, we perform the Morawetz type energy estimates,
utilizing the Fefferman-Stein inequality on the duality between ${\cal H}^1$
and the space of functions of bounded mean oscillation.
Another diffusion wave part $M(t)-{\cal K}_{\nu} (t) \ast M_{0,in}$
is shown to grow at the rate of order $\log t$ as $t$ goes to infinity.
Here we use  fundamental solutions for the linearized system given in \cite{Kobayashi-Tsuda}.  Since a high frequency part of the solutions
exponentially decays,
a low frequency part has only to be estimated here.
By direct computation with the explicit form
of the Green matrix
we get the growth order for $M(t)-{\cal K}_{\nu} (t) \ast M_{0,in}$.
For the Stokes flow part ${\cal K}_{\nu}(t)\ast M_{0,in}$,
the space-time $L^2$ boundedness is derived
in Theorem \ref{hth} bellow. 
Combining these estimates for diffusion wave
and the Stokes flow parts yields
the asymptotic expansion \eqref{behavior}.

  This paper is organized as follows. In section 2 some notations and lemmas are given.
In section 3, the main results  are presented.
In section 4, the proofs of the estimates for the diffusion wave parts are demonstrated.

\section{Preliminaries}
In this section we introduce some notations such as function spaces,
used in this paper.
We also present some lemmas,
needed in the proof of the main result.

The norm on $X$ is denoted by $\|\cdot\|_{X}$ for a given Banach space $X$.

Let $1\leqq p \leqq \infty.$
$L^p$ is
the usual Lebesgue
space of $p-$th powered integrable and essentially bounded
functions on $\mathbb{R}^2$
for a finite $p$ and $p = \infty$, respectively.
Let $k$ be a nonnegative integer.
$W^{k,p}$ and $H^k$ are the usual
Sobolev space of order $k$,
based on $L^2$ and $L^p$, respectively.  
As usual, $H^0$ is defined by $H^0:=L^2$. 

We also use the notation $L^p$ to denote the function space 
of all vector fields 
$w=\trans(w_1, w_2)$ on $\mathbb{R}^2$
satisfying
$w_{j}\in L^p$ $(j=1, 2)$,
and $\|\cdot\|_{L^p}$ is
the norm $\|\cdot\|_{(L^p)^2}$ for brevity
if  no confusion will occur.
Similarly, a function space $X$ is
the linear space
of all vector fields 
$w=\trans(w_1, w_2)$ on $\mathbb{R}^2$
satisfying
$w_{j}\in X$ $(j=1, 2)$,
and  $\|\cdot\|_{X}$ is the norm $\|\cdot\|_{X^2}$ if no confusion will occur. 

Let $u=\trans(\phi,m)$ with $\phi\in H^k$ and $m=\trans(m_1, m_2)\in H^j$.  
Then the norm $\|u\|_{H^k\times H^j}$
is defined as that of $u$ on $H^k\times H^j$
$$
\|u\|_{H^k\times H^j}:=\left(\|\phi\|_{H^k}^2+\|w\|_{H^j}^2\right)^{\frac{1}{2}}.
$$
In particular, if $j=k$,  we put
$$
H^k:=H^k\times (H^k)^2, 
\ \ \ 
\|u\|_{H^k}:=\|u\|_{H^k\times (H^k)^2} 
\ \ \ (u=\trans(\phi,m)).
$$
Let X and Y be given Banach spaces. For $u=\trans(\phi,m)\in X\times Y$ with $m=\trans(m_1, m_2)$, similarly we set 
$$
\|u\|_{X\times Y}:=\left(\|\phi\|_{X}^2+\|m\|_{Y}^2\right)^{\frac{1}{2}}.
$$
More generally, in the case that $Y=X^2$, let
$$
X:=X\times X^2, 
\ \ \ 
\|u\|_{X}:=\|u\|_{X\times X^2}.
$$


The symbols $\hat{f}$ and $\mathcal{F}[f]$ stand for  the Fourier transform of $f$
with respect to
the space variables $x$
\begin{eqnarray}
\hat{f}(\xi)
=\mathcal{F}[f](\xi)
:=\int_{\mathbb{R}^2}f(x)e^{-ix\cdot\xi}dx,
\quad \xi\in\mathbb{R}^2.
\nonumber
\end{eqnarray}
Furthermore,  the inverse Fourier transform of $f$ is defined by
\begin{eqnarray}
\mathcal{F}^{- 1}[f](x)
:=(2\pi)^{-2}\int_{\mathbb{R}^2}f(\xi)e^{i\xi\cdot x}d\xi,
\quad
x\in\mathbb{R}^2.
\nonumber
\end{eqnarray}

%
%
%
%
%
%
%
%
For a nonnegative number $s$, $[s]$ is the Gaussian symbol which denotes the integer part of $s$. 
The symbol $``\ast"$ denotes the convolution on the space variable $x$. 
\par
\vskip5pt
\parindent0pt
$\quad$
Now we prepare the Hardy space ${\cal H}^1$ and BMO space. 

\vspace{2ex}

\begin{defi}
The Hardy space ${\cal H}^1 = {\cal H}^1  (\mathbb{R}^2)$ consists of integrable functions on $\mathbb{R}^2$,
	$f \in L^1(\mathbb{R}^2)$ such that
\begin{eqnarray}
		\|f\|_{{\cal H}^1(\mathbb{R}^2)}
		= \int_{\mathbb{R}^{2}}
		\sup \limits_{r \textgreater 0}
		|\phi_r \ast f(x)| dx
		\nonumber
	\end{eqnarray}
	is finite, where
	$\phi_r(x) = r^{-n}\phi(r^{-1}x)$
	for $r \textgreater 0$ and $\phi$ is
	a smooth function on $\mathbb{R}^2$ with
	compact support in an unit ball with
	center of the origin,
	$B_1(0)=\{x \in \mathbb{R}^2 ; |x|\textless 1\}$. 
The definition dose not depend on choice of a function $\phi$.

\end{defi}

\vspace{2ex}

\begin{defi} 
Let $f$ be a locally integrable in $\mathbb{R}^2$,
 $f \in L^{1}_{loc} (\mathbb{R}^2)$.
	We say that $f$ is of bounded mean oscillation,
	abbreviated as $BMO$, if
	\begin{eqnarray}
		\|f\|_{{\rm BMO}}
		= \sup \limits_{B \subset \mathbb{R}^2}
		\frac{1}{|B|}
		\int_B
		|f - (f)_B|
		dx
		\textless \infty,
		\nonumber
	\end{eqnarray}
	where the supremum ranges over all finite ball
	$B$ $\subset$ $\mathbb{R}^2$,
	$|B|$ is the $2$-dimensional Lebesgue measure of $B$,
	and $(f)_B$ denotes the integral mean of $f$ over $B$,
	namely
	$(f)_B = \frac{1}{|B|} \ \int_B f(x) dx$.
	\par
	\vskip5pt
	\parindent0pt
	$\quad$
	The class of functions of $BMO$, modulo constants,
	is a Banach space with the norm
	$\| \cdot \|_{BMO}$
	defined above.
	
\end{defi}

\vspace{2ex}

We crucially use
the decisive Fefferman-Stein inequality,
which means the duality between
${\cal H}^1(\mathbb{R}^2)$ and $BMO(\mathbb{R}^2)$, i.e., 
$({\cal H}^1(\mathbb{R}^2))^{\ast} = BMO(\mathbb{R}^2)$.
For the proof, see \cite{Fefferman-Stein}.

\vspace{2ex}

\begin{lem}
\label{Fefferman-Stein inequality}
{\rm (Fefferman-Stein inequality)}
	There is a positive constant $C$
	 such that if
	$f \in {\cal H}^1(\mathbb{R}^2)$
	and
	$g \in BMO(\mathbb{R}^2)$,
	then
	\begin{eqnarray}
		\left|
		\int_{\mathbb{R}^2}
		f g
		dx
		\right|
		\le
		C
		\|f\|_{{\cal H}^1 (\mathbb{R}^2)}
		\|g\|_{BMO (\mathbb{R}^2)}
		\nonumber.
	\end{eqnarray}
\end{lem}

\vspace{2ex}

$\quad$
We also recall
the well known Poincar\'{e} inequality.

\vspace{2ex}

\begin{lem}
\label{poancare}
It holds that 
$$
\|f\|_{BMO(\mathbb{R}^2)} \leq \|\nabla f \|_{L^2 (\mathbb{R}^2)}
$$
for $f \in H^1 (\mathbb{R}^2)$.
\end{lem}

\vspace{2ex}

$\quad$We denote by $C^{\infty}_{0, \sigma}$ the set of all vector valued functions $\phi=\trans(\phi_1, \phi_2)$ whose each $\phi_j$ $(j=1,2)$  is $C^\infty$ function having compact support, and satisfying that $\div \phi=0$.  For $1\leq q<\infty$,  $L^q_\sigma$ is the closure of  $C^{\infty}_{0, \sigma}$ with respect to the $L^q$ norm.  

$\quad$A spatial weighted function space
$W^{1,2}_w (\R^2)$ is defined by 
\begin{align*}
W^{1,2}_w (\R^2)=
\left\{
u\, : \, \frac{u}{w(x)}\in L^2(\R^2),
\nabla u\in L^2(\R^2)
\right\},
\end{align*}
where $w(x)$ is a spatial weight defined by $w(x)=
(1+|x|)\log (2+|x|)$.

The following H\"{o}lder type inequality is proved by Amrouche and Nguyen \cite{AN11}. 

\vspace{2ex}

\begin{lem} \label{1AN}
{\rm (\cite[Corollary 2.10]{AN11})}
Let $f\in L^1_{\sigma}(\R^2)$. Then it holds true that $\int_{\R^2}f(x)dx=0$ and that, for such $f$ and any $g\in W^{1,2}_w (\R^2) \cap L^\infty (\R^2)$, 
\begin{align*}
\left| 
\int_{\R^2}fgdx 
\right| 
\leq
C\|f\|_{L^1}\|\nabla g\|_{L^2}. 
\end{align*}
\end{lem}
%
Since $W^{1, 2} (\R^2) \subset W^{1, 2}_w (\R^2)$,  Lemma \ref{1AN} also yields the following  
\begin{cor} \label{AN}
Let $f\in L^1_{\sigma}(\R^2)$. Then there holds that $\int_{\R^2}f(x)dx=0$ and, for such $f$ and any $g\in W^{1,2} (\R^2) \cap L^\infty (\R^2)$, 
\begin{align*}
\left| 
\int_{\R^2}fgdx 
\right| 
\leq
C\|f\|_{L^1}\|\nabla g\|_{L^2}. 
\end{align*}

\end{cor}
 
\section{Main results}
In this section, we consider the linearized system
corresponding to (\ref{CNSK})
and present some decay estimates for its solution.
One of key estimates to show \eqref{behavior} is a space-time $L^2$ boundedness of
the density for the linearized system.
First of all, 
\eqref{CNSK} is reformulated and linearized as follows.  
Hereafter we assume that $\rho_*=1$ without loss of generality.
We also set 
$$
\phi=\rho-1,
\quad
m=\frac{M}{\gamma},
\quad
\gamma=\sqrt{P'(1)}.
$$
Substituting $\phi$ and $m$ into \eqref{CNSK},
we have the system of equations
\begin{eqnarray}
\left\{
\begin{array}{lll}
\partial_{t}\phi +\gamma\div  m=0,
\\
\partial_{t}m-\nu\Delta m-\tilde{\nu}\nabla\div m+\gamma \nabla \phi-\kappa_0 \nabla \Delta\phi=f(u),
\\
\phi|_{t=0}=\phi_0,
\ \ m|_{t=0}=m_0,
\label{cnsk-nolinear}
\end{array}
\right.
\end{eqnarray}
where we use the notation
$$
u=\trans(\phi,m),
\quad
\phi_0={\rho_0-1},
\quad
m_0=\frac{M_0}{\gamma},
\quad
\nu={\mu},
\quad
\tilde{\nu}={\mu+\mu'},
\quad
\kappa_0=\frac{\kappa}{\gamma}
$$
and put 
\begin{eqnarray*}
f(u)&=& -\Big\{\gamma\div(m \otimes m)+\gamma\div(P_{(1)}(\phi)\phi m \otimes m)+\frac{1}{\gamma}\nabla(P_{(2)}(\phi)\phi^2)\\
&&\qquad -\nu\Delta(P_{(1)}(\phi)\phi m)-\tilde{\nu}\nabla\div(P_{(1)}(\phi)\phi m)-\div \Phi(\phi)\Big\},\\
P_{(1)}(\phi)
&=&
\int_{0}^{1}f'(1+\tau\phi)d\tau,
\quad
f(\tau)=\frac{1}{\tau},
\quad \tau \in \mathbb{R},
\\
P_{(2)}(\phi)
&=&
\int_{0}^{1}(1-\tau)P''\Big(1+\tau\phi\Big)d\tau,\\
\Phi(\phi)
&=&\kappa_0\Big\{\phi\Delta \phi I_2 +(\nabla \phi)\cdot(\nabla \phi)I_2 -\frac{|\nabla \phi |^2}{2}I_2-\nabla \phi \otimes \nabla \phi \Big\}.
\end{eqnarray*}
Therefore, $(\ref{CNSK})$ is linearized as 
\begin{eqnarray}
\left\{
\begin{array}{lll}
\partial_{t}\phi +\gamma\div  m=0,
\\
\partial_{t}m-\nu\Delta m-\tilde{\nu}\nabla\div m+\gamma \nabla \phi-\kappa_0 \nabla \Delta\phi=0,
\\
\phi|_{t=0}=\phi_0,
\ \ m|_{t=0}=m_0.
\label{cnsk-linear}
\end{array}
\right.
\end{eqnarray}
By \eqref{cnsk-linear}, $\phi$ satisfies the following doubly dissipative equation 
\begin{numcases}
{}
\partial_{tt} \phi -(\nu+\tilde{\nu})\Delta \phi_t - \gamma^2 \Delta \phi + \gamma \kappa_0 \Delta^2 \phi=0,
\label{visco+}\\
\phi(x,0) =\phi_0, \ \ \del_t \phi(x,0)=-\gamma \div m_0.
\nonumber
\end{numcases}
Due to the positivity of $\nu$ and $\tilde{\nu}$,
we may suppose that $\nu+\tilde{\nu}=1$ and $\gamma=1$
without loss of generality. Then $\phi$ satisfies 
\begin{numcases}
{}
\partial_{tt} \phi -\Delta \phi_t - \Delta \phi + \kappa_0 \Delta^2 \phi=0,
\label{visco+-normalize}\\
\phi(x,0) =\phi_0, \ \ \del_t \phi(x,0)=-\div m_0.
\nonumber
\end{numcases}

\vspace{2ex}
%
%
%
%
\par
\parindent0pt
$\quad$
Now we state the existence of solutions to \eqref{visco+-normalize}
in the energy class, defined in the following. 

\vspace{2ex}

\begin{defi}
A function $\phi$ defined on $(0, \infty) \times \mathbb{R}^2$ is called to be a
solution to \eqref{visco+-normalize} 
if $\phi$ belongs to $L^\infty(0,\infty; H^2) $ with $\del_t \phi \in L^\infty(0,\infty; L^2)$ and satisfies \eqref{visco+-normalize} in the distribution sense.
\end{defi}

\vspace{2ex}

\begin{thm}\label{weak-sol1}
For each $(\phi_0, m_0) \in H^2 \times H^1$ there exists a unique solution $\phi \in L^\infty(0,\infty; H^2) $ with $\del_t \phi \in L^\infty(0,\infty; L^2)$ to \eqref{visco+-normalize} such that  
\begin{eqnarray} 
 &&\frac{1}{2}(\|\phi_t(t)\|_{L^2}^2+\|\nabla \phi(t)\|_{L^2}^2+\kappa_0\|\Delta\phi(t)\|_{L^2}^2)+ 
  \displaystyle\int_0^t\|\nabla \phi_\tau(\tau)\|_{L^2}^2d\tau \nonumber\\
 &&=  \frac{1}{2}(\|\div m_0\|_{L^2}^2+\|\nabla \phi_0\|_{L^2}^2+\kappa_0\|\Delta\phi_0\|_{L^2}^2)\label{ene-weak-sol}
 \end{eqnarray}
holds for any $t \geq 0$. 
\end{thm}

\vspace{2ex}
Theorem \ref{weak-sol1} is valid by the standard Galerkin method based on the energy inequality \eqref{ene-weak-sol} in a similar manner to  the proof of Lemma 2.3 in Huafei and Yadong \cite{Huafei-Yadong} and we omit the details. 
 \vspace{2ex}

$\quad$First, we show a boundedness
in $L^2(0,\infty ; L^2)$
for a solution $\phi$ to \eqref{visco+-normalize}.

\vspace{2ex}

\begin{thm}\label{time-spatial est}
Suppose that $\phi_0 \in H^2 \cap {\cal H}^1$, $m_0 \in H^1$ and $m_0+ \nabla \phi_0 \in  {\cal H}^1$. Set 
$$
J_0 = (\kappa_0+\kappa_0^2)\{\|\trans(\phi_0,m_0)\|_{H^1}^2+\|\Delta \phi_0\|_{L^2}^2\}+(1+\kappa_0)\|m_0+\nabla \phi_0\|_{ {\cal H}^1}^2 +\|\phi_0\|_{{\cal H}^1\cap L^2}^2. 
$$
Let $\phi$ be a solution to $(\ref{visco+-normalize})$.
Then it holds true that
$$
\displaystyle\int_{0}^{t}\|\phi(\tau)\|_{L^2}^2 d\tau \leq CJ_0
$$
for any $t > 0$, where $C$ is a positive constant independent of $t$ and $\kappa_0$. 
\end{thm}

\vspace{2ex}

$\quad$
We can also treat the linearized CNS, that is (\ref{visco+-normalize}) with
the zero capillary constant, $\kappa_0 = 0$.
\begin{numcases}
{}
\partial_{tt} \phi -\Delta \phi_t - \Delta \phi = 0,
\label{2visco+-normalize}\\
\phi(x,0) =\phi_0, \ \ \del_t \phi(x,0)=-\div m_0.
\nonumber
\end{numcases}

\vspace{2ex}

\begin{defi}
A function $\phi$ defined on $(0, \infty) \times \mathbb{R}^2$ is called to be a
solution to \eqref{2visco+-normalize}
if $\phi$ belongs to $C([0,\infty); H^1) \cap C^1([0,\infty); L^2) $ and satisfies \eqref{2visco+-normalize} in the distribution sense.
\end{defi}

\vspace{2ex}

The existence of a unique solution to \eqref{2visco+-normalize} is well known as follows.
For the proof we can refer to Proposition 2.1 in Ikehata, Todorova and Yordanov \cite{Ikehata} using the Lumer-Phillips theorem.

\vspace{2ex}

\begin{thm}\label{weak-sol2}
For each $(\phi_0, m_0) \in H^1$ there exists a unique solution $\phi \in C([0,\infty); H^1) \cap C^1([0,\infty); L^2) $ to \eqref{2visco+-normalize} such that  
\begin{eqnarray*} 
 &&\frac{1}{2}(\|\phi_t(t)\|_{L^2}^2+\|\nabla \phi(t)\|_{L^2}^2+
  \displaystyle\int_0^t\|\nabla \phi_\tau(\tau)\|_{L^2}^2d\tau \\
 &&=  \frac{1}{2}(\|\div m_0\|_{L^2}^2+\|\nabla \phi_0\|_{L^2}^2)\label{2ene-weak-sol}
 \end{eqnarray*}
holds for any $t \geq 0$. 
\end{thm}

 \vspace{2ex}

%
%
%
%
In the case that $\kappa_0 = 0$, 
we also have the time-space $L^2$ estimate for linearized CNS.


\vspace{2ex}

\begin{thm}\label{time-spatial est-CNS}
Let $(\phi_0, m_0) \in H^1$, $\phi_0 \in {\cal H}^1$ and $m_0+ \nabla \phi_0 \in  {\cal H}^1$. Set 
$$
J_1 = \|m_0+\nabla \phi_0\|_{ {\cal H}^1}^2 +\|\phi_0\|_{{\cal H}^1\cap L^2}^2. 
$$
Let $\phi$ be a solution to $(\ref{2visco+-normalize})$.
Then there holds
$$
\displaystyle\int_{0}^{t}\|\phi(\tau)\|_{L^2}^2d\tau \leq C J_1
$$
for any $t > 0$, where $C$ is a positive constant independent of $t$.  
\end{thm}

\vspace{2ex}

\par
\parindent0pt
$\quad$
Next, we have a time decay estimate of the solution in the energy class to \eqref{visco+-normalize}. Note that by Theorem \ref{weak-sol1} and the Sobolev inequality $\phi \in C([0,\infty); L^2)$. We have the following

\vspace{2ex}

\begin{thm}\label{time-decay est-CNSK}
Under the assumption of Theorem \ref{time-spatial est}, it holds that
$$
(1+t)\|\phi(t)\|_{L^2}^2  \leq C J_0
$$
for any $t > 0$, where $C$ is a positive constant independent of $t$.  
\end{thm}

\vspace{2ex}
%
%
%
%

\par
\parindent0pt
$\quad$
We
now recall the existence of solutions to linear system \eqref{cnsk-linear} in the energy class in order to consider another diffusion wave part
$m-{\cal K}_{\nu}\ast m_{0,in}$.  
The system \eqref{cnsk-linear} is rewritten as 
\eqn{ns2}
$$
\delt u+Au=0,
$$
where 
\begin{eqnarray}\label{(a)}
u=\trans(\phi, m), \ \ A=\begin{pmatrix}
0 &\gamma\mathrm{div}\\
\gamma\nabla -\kappa_0 \nabla \Delta &-\nu\Delta-\tilde{\nu}\nabla\mathrm{div}\\
\end{pmatrix}
\end{eqnarray}
Let  us introduce a semigroup $S(t)=e^{-tA}$ generated by $A$; 
  $$
  S(t)=e^{-tA}={\cal F}^{-1}e^{-t\hat A_{\xi}}{\cal F},
  $$
where 
\begin{eqnarray}
\hat{A}_{\xi}=\begin{pmatrix}
0 &i\gamma\trans{\xi}\\
i\gamma \xi+i\kappa_0|\xi|^2 \xi &\nu|\xi|^{2}I_{n}+\tilde{\nu}\xi \trans{\xi}
\end{pmatrix}
\ \ \ (\xi\in \mathbb{R}^2). 
\nonumber
\end{eqnarray}

\vspace{2ex}

\begin{thm}{\rm \label{existence-sol-linear}\cite[Proposition 3.3]{Tsuda-cnsk}
Let $s$ be a nonnegative integer satisfying $s\geq2$. Then 
$S(t)=e^{-tA}$ is a contraction semigroup for \eqref{ns2} on $H^s\times H^{s-1}$. 
In addition, for each
$u_0 = \trans(\phi_0, m_0) \in H^s\times H^{s-1}$
and all
$T > 0$,
$S(t)$ satisfies 
$$
S (\cdot)u_0 \in C([0,T];H^s\times H^{s-1}), \ S(0)u_0=u_0
$$
and there holds the estimate
\begin{eqnarray}
\| S (t)u_0\|_{H^s\times H^{s-1}}\leq 
\|u_0\|_{H^s\times H^{s-1}}\label{estimate-S(t)}
\end{eqnarray}
for $u_0 =  \trans(\phi_0, m_0) \in H^s\times H^{s-1}$ and $t\geq 0$. }
\end{thm}

\vspace{2ex}

\begin{rem}{\rm 
Proposition 3.3 in \cite{Tsuda-cnsk} is stated on the three dimensional case. However,
the proof is based on the standard energy estimate for the resolvent problem in the Fourier space
and it can be also applied to
our two dimensional case.  }
\end{rem}

\vspace{2ex}

Finally, another diffusion wave part
$m-{\cal K}_{\nu}\ast m_{0,in}$ is shown to grow
in $L^2(0,\infty; L^2)$
at the rate of logarithmic order.

\vspace{2ex}

%
%

\begin{thm}\label{time-spatial est-diffusive2}
Let $u_0 = \trans(\phi_0, m_0) \in H^2\times H^1$
and $u$ be a solution of \eqref{cnsk-linear},
$u (t) = S (t) u_0$, as in Theorem \ref{existence-sol-linear}.
Suppose that $m_0 \in L^2 \cap L^1$, $|x| m_0 \in L^1$ and $\hat{m}_0 (0) \neq 0$.
Then it holds true that 
$$
\int_{0}^{t}\|m(\tau)-{\cal K}_{\nu}(\tau)\ast m_{0,in}\|_{L^2}^2d\tau \sim \log t  \ \ \mbox{as}  \ \ t \rightarrow \infty, 
$$
precisely,
$$
\limsup_{t\rightarrow \infty}\Big|\displaystyle\int_{0}^{t}\|m(\tau)-{\cal K}_{\nu}(\tau)\ast m_{0,in}\|_{L^2}^2d\tau -\log t \Big|\leq C,
$$
where $C=C(u_0)$ is a positive constant independent of $t$. 
\end{thm}

\vspace{2ex}

\begin{rem} {\rm 
In addition of the initial condition
in Theorem \ref{time-spatial est},
we assume that $m_0 \in {\cal H}^1$ and then, it holds that 
\begin{eqnarray}
\displaystyle
\int_{0}^{t}\|m(\tau)-{\cal K}_{\nu}(\tau)\ast m_{0,in}\|_{L^2}^2 d\tau
\leq C,
\quad t > 0,
\label{strong-assump}
\end{eqnarray}
where $C=C(u_0)$ is  a positive constant independent of $t$. 
This shows a gain of regularity by the membership in Hardy space of data,
similarly as in the decay estimates for density in Theorems \ref{time-spatial est}
and \ref{time-spatial est-CNS}.
The similar phenomenon have already observed in
\cite{Kobayashi-Misawa,Kobayashi-Misawa2}
for the dissipative wave equations.
The proof is given by direct computations based on the explicit form of fundamental solution
\eqref{convolution-solutionfomula}
below and a similar argument
as in Kobayashi and Misawa \cite{Kobayashi-Misawa, Kobayashi-Misawa2}.
We omit the details, here.
}
\end{rem}

\vspace{2ex}

\quad We shall state the space-time $L^2$ boundedness for the Stokes flow part ${\cal K}_{\nu}(t)\ast M_{0,in}$.   
In order to state the boundedness let us introduce the following incompressible Stokes system 
\begin{align}
\begin{cases}\label{SE}
\pt u-\nu \Delta u +\nabla P
=0
&\quad {\rm in}\ (0,\infty )\times \R^2,\\
\div u=0
&\quad {\rm in}\ (0,\infty )\times \R^2,\\
u(0,x)=M_{0,in}
&\quad {\rm in}\ \R^2.
 \end{cases}
\end{align}
The following Helmholtz decomposition is well known. (Cf., Simader and Sohr \cite{SS}) 
$$
L^q (\R^2)= L^q_\sigma (\R^2) \oplus G^q (\R^2), \ \ (1<q< \infty), 
$$ 
where $G^q (\R^2)$ denotes the set of all functions of the potential flow part and defined by $G^q=\{\nabla p \in L^q (\R^2); p\in L^q_{{\rm loc}}(\R^2)\}$.  Here we denote by $P_q$ the projection operator from $L^q (\R^2)$ to  $L^q_\sigma (\R^2)$. 
 On the whole space $P_q$ is given by the Riesz operator 
$$
P_q f =\mathcal{F}^{-1}\Big\{\Big(I_2-\frac{\xi\trans{\xi}}{|\xi|^2}\Big)\hat{f}\Big\}. 
$$
Applying  the Helmholtz projection to the 
Stokes  equations \eqref{SE} derives the following system 
\begin{align}
\begin{cases}\label{SE2}
\pt u-\nu P_q \Delta u
=0
&{\rm in}\ (0,\infty )\times \R^2,\\
\div u=0
&{\rm in}\ (0,\infty )\times \R^2,\\
u(0,x)=u_0:=M_{0,in}
&{\rm in}\ \R^2.
\end{cases}
\end{align}

We define the Stokes operator $A_q$ on $L^q_{\sigma}$ by $A_q =-P_q \Delta$ with domain $D(A_q)= W^{2,q}(\R^2)\cap L^q_\sigma (\R^2)$.  Concerning existence of solutions to \eqref{SE2}, we heve

\vspace{2ex}
\begin{thm}{\rm (Giga and Sohr \cite{Giga-Sohr})}  $-A_q$ generates a uniformly bounded holomorphic semigroup $\{e^{-t A_q}\}_{t\geq 0}$ of class $C^0$ in  $L^q_\sigma (\R^2)$. 
\end{thm}
\vspace{2ex}

Note that  the solution to \eqref{SE2} satisfies $u(t)={\cal K}_{\nu}(t)\ast M_{0,in}$ and we can thus estimate the solution $u$ to \eqref{SE2} as follows 

\begin{thm} \label{hth}
Let $M_{0,in}\in L^2 \cap L^1$ and $u$ be a solution to \eqref{SE2}. Then $u$ satisfies the estimate 
\begin{align*}
\displaystyle\int_0^t\|u(s)\|_{L^2}^2ds\leq C\|M_{0,in}\|_{L^1}^2
\end{align*}
uniformly for $t$.  
\end{thm}

\vspace{2ex}

We give the proof of Theorem \ref{hth} here. 
\vspace{2ex}

\textbf{Proof of Theorem \ref{hth}.} Put $v(t,x):=\int_0^t u(s,x)ds.$  Then $v$ satisfies

\begin{align}
\begin{cases}
\displaystyle
\pt v-\nu \Delta v=u_0&{\rm in}\ (0,\infty )\times \R^2,\\\
v(0,x)=0 &{\rm in}\ \R^2.
\label{3.1}
\end{cases}
\end{align}
Here we used that $P_2 (\Delta u)=\Delta (P_2 u)=u$ in \eqref{SE2}.  
A test function $\pt w$ in \eqref{3.1}, being integrated on the time interval $(0, t)$, and the fact $\del_t w=u$ yield the estimate 

\begin{equation}
\int_0^t\|u(s)\|_{L^2}^2ds
+
\frac{\nu}{2}\|\nabla v(t)\|_{L^2}^2
=
(u_0,v(t))
\label{3.2}
\end{equation}
The first term of the right hand
side in \eqref{3.2} is estimated by  Corollary \ref{AN} as follows
\begin{align}
|(u_0,v(t))|
\leq 
C\|u_0\|_{L^1}\| \nabla v(t)\|_{L^2}&\leq
C_1\|u_0\|_{L^1}^2
+\frac{\nu}{4}\| \nabla v(t)\|_{L^2}^2. 
\label{PL1}
\end{align}
\eqref{3.2} and \eqref{PL1} derive the desired estimate.  
The proof is completed. $\hfill\square$

\vspace{2ex}

From \eqref{strong-assump} together with
Theorem \ref{hth} and \cite[Chapter 3, Section 3, Theorem 3]{Stein} 
we find that
if
$m_0 \in {\cal H}^1$ is
added in
the assumption of Theorem \ref{time-spatial est},
the space-time $L^2$ boundedness holds true for  $M(t)-{\cal K}_{\nu}\ast M_{0,in}(t)$, $\rho(t)- \rho_*$ and ${\cal K}_{\nu}(t)\ast M_{0,in}$.

\section{Proof of main results}

\subsection{Proof of  Theorems \ref{time-spatial est} and \ref{time-decay est-CNSK}} 

In this subsection, we prove Theorems \ref{time-spatial est} and \ref{time-decay est-CNSK}.
The proof is performed
by modifying Morawetz's energy estimate.
For a solution $\phi$ to \eqref{visco+-normalize}, we define the function $w$ by 
\begin{eqnarray}
w = \displaystyle\int_{0}^{t} \phi(\tau) d\tau -\div \Phi, \label{w-okikata}
\end{eqnarray}
where $\Phi$ is a solution to the doubly Laplace equation 
\begin{eqnarray}
(-\Delta + \kappa_0 \Delta^2) \Phi =m_0 + \nabla \phi_0.
\label{doubly-eqn}
\end{eqnarray}
\par
\parindent0pt
$\quad$
For the existence of solution to (\ref{doubly-eqn}) we have
\begin{thm}
\label{exist-d-eqn}
Suppose that
$\phi_0 \in H^2 \cap {\cal H}^1 (\mathbb{R}^2)$,
$m_0 \in H^1 (\mathbb{R}^2)$
and $m_0 + \nabla \phi_0 \in {\cal H}^1 (\mathbb{R}^2)$.
Then, there exists a solution $\Phi$ of
\begin{eqnarray}
(-\Delta) \Phi =(1-\kappa_0 \Delta)^{-1}(m_0 + \nabla \phi_0)
\quad \mbox{in $\mathbb{R}^2$}
\label{2-d-eqn}
\end{eqnarray}
such that
\begin{eqnarray}
&&
\|\nabla \Phi\|_{L^2 (\mathbb{R}^2)}
+
\|\Delta \Phi\|_{L^2 (\mathbb{R}^2)}
\leq C\|m_0 + \nabla \phi_0\|_{{\cal H}^1 (\mathbb{R}^2)},
\label{Phi-est1}
\\
&&
\|\div \Delta \Phi\|_{L^2 (\mathbb{R}^2)}
\le
C(\|\Delta \phi_0\|_{L^2 (\mathbb{R}^2)}+\|\trans(\phi_0, m_0)\|_{H^1 (\mathbb{R}^2)}).
\label{Phi-est2}
\end{eqnarray}
\end{thm}
\par
\parindent0pt
The proof of Theorem \ref{exist-d-eqn} is in Appendix.
\par
\parindent0pt
$\quad$
By definition of $w$ we derive 
\begin{eqnarray}
\left\{
\begin{array}{ll}
w_{tt}-\Delta w_t -\Delta w +\kappa_0 \Delta^2 w=0,
\\
w(0)=-\div \Phi, \ \ w_t (0)=\phi_0.
\end{array}
\right.\label{w}
\end{eqnarray}
To estimate 
$$
\displaystyle\int_0^t \|\phi(\tau)\|_{L^2}^2d\tau = \displaystyle\int_0^t \|w_\tau(\tau)\|_{L^2}^2 d\tau,
$$
we take $L^2$ inner product of $(\ref{w})_1$ with $w$ and thus, have
\begin{eqnarray*}
\frac{d}{dt}(w_t, w)+\frac{1}{2}\frac{d}{dt}\|\nabla w\|_{L^2}^2+\|\nabla w\|_{L^2}^2+\kappa_0 \|\Delta w\|_{L^2}^2=\|w_t\|_{L^2}^2,
\end{eqnarray*}
which is integrated on the time interval $(0, t)$, yielding
\begin{eqnarray}
\lefteqn{\displaystyle\int_0^t \|w_\tau(\tau)\|_{L^2}^2 d\tau}\nonumber\\
&&= \displaystyle\int_0^t\|\nabla w(\tau)\|_{L^2}^2 d\tau +\displaystyle\int_0^t \kappa_0\|\Delta w(\tau)\|_{L^2}^2 d\tau\nonumber\\
&&\qquad +\frac{\|\nabla w(t)\|_{L^2}^2}{2}+(w_t, w)  - \frac{\|\nabla w (0)\|_{L^2}^2}{2}-(w_t(0), w(0)),\label{w-base} 
\end{eqnarray}
where we note that $w(0)=-\div \Phi$. 
Now we will estimate the terms in the right hand side of (\ref{w-base}). The terms $\|\nabla w(0)\|_{L^2}^2$ and $(w_t (0), w(0))$ are directly estimated by \eqref{Phi-est1}. 
A test function $w_t$ in \eqref{w}, being integrated on the time interval $(0, t)$
and using \eqref{Phi-est1} and \eqref{Phi-est2}, yield the estimate 
 \begin{eqnarray}
\lefteqn{\|w_t(t)\|_{L^2}^2+\|\nabla w(t)\|_{L^2}^2+\kappa_0\|\Delta w (t)\|_{L^2}^2 +\displaystyle\int_0^t
\|\nabla w_\tau(\tau)\|_{L^2}^2 d\tau}\nonumber\\
&&\leq C(\|\phi_0\|_{L^2}^2 +\|m_0+\nabla \phi_0\|_{{\cal H}^1}^2)+C\kappa_0(\|\Delta \phi_0\|_{L^2}^2+\|\trans( \phi_0, m_0)\|_{H^1}^2). \label{nabla-w-est}
\end{eqnarray}
 
 On the other hand,
from taking $L^2$ inner product of $-\kappa_0 \Delta w$ with \eqref{w}
and integrating on $(0,t)$ we obtain that 
\begin{eqnarray*}
&&\frac{1}{2}\kappa_0\|\Delta w(t)\|_{L^2}^2+\kappa_0\displaystyle\int_0^t\|\Delta w(\tau)\|_{L^2}^2d\tau 
+\kappa_0(\nabla w_t, \nabla w)+ \kappa_0^2\displaystyle\int_0^t\|\nabla \Delta w(\tau)\|_{L^2}^2d\tau\\
&&\quad = \frac{1}{2}\kappa_0 \|\Delta w(0)\|_{L^2}^2+\kappa_0 \displaystyle\int_0^t\|\nabla w_\tau(\tau)\|_{L^2}^2d\tau +\kappa_0(\nabla w_t(0), \nabla w(0))
\end{eqnarray*}
and thus,
\begin{eqnarray} 
\lefteqn{\frac{1}{2}\kappa_0\|\Delta w(t)\|_{L^2}^2+\kappa_0\displaystyle\int_0^t\|\Delta w(\tau)\|_{L^2}^2d\tau 
+ \kappa_0^2\displaystyle\int_0^t\|\nabla \Delta w(\tau)\|_{L^2}^2d\tau}\nonumber\\
&&\leq C\kappa_0\{\|\nabla \phi_0\|_{L^2}^2+\|\nabla \div \Phi\|_{L^2}^2+\|\Delta \div \Phi\|_{L^2}^2+
\displaystyle\int_0^t
\|\nabla w_\tau(\tau)\|_{L^2}^2 d\tau\}\nonumber\\
&&\qquad +C\kappa_0(\|\nabla w(t)\|_{L^2}^2+\|\nabla w_t(t)\|_{L^2}^2). \label{Delta-w-est}
 \end{eqnarray}
The terms having $\Phi$ in \eqref{Delta-w-est} are estimated by
\eqref{Phi-est1} and \eqref{Phi-est2}.
The 4th and 5th terms in the right hand side of (\ref{Delta-w-est})
are also estimated by \eqref{nabla-w-est}. 
For the term $\|\nabla w_t(t)\|_{L^2}^2=\|\nabla \phi(t)\|_{L^2}^2$,
we apply the standard energy estimate, obtained from \eqref{visco+-normalize}
with a test function $\phi_t$
  \begin{eqnarray*} 
 \lefteqn{\frac{1}{2}(\|\phi_t(t)\|_{L^2}^2+\|\nabla \phi(t)\|_{L^2}^2+\kappa_0\|\Delta\phi(t)\|_{L^2}^2)+ 
  \displaystyle\int_0^t\|\nabla \phi_\tau(\tau)\|_{L^2}^2d\tau }\\
 &&=  \frac{1}{2}(\|\phi_t(0)\|_{L^2}^2+\|\nabla \phi(0)\|_{L^2}^2+\kappa_0\|\Delta\phi(0)\|_{L^2}^2)\\
 &&\leq  C(\|\div m_0\|_{L^2}^2+\|\nabla \phi_0\|_{L^2}^2)+C\kappa_0 \|\Delta \phi_0\|_{L^2}^2
 \end{eqnarray*}
and thus,
 \begin{eqnarray} 
 \kappa_0\displaystyle\int_0^t\|\Delta w(\tau)\|_{L^2}^2d\tau &\leq& C(\kappa_0+\kappa_0^2) \{\|\trans(\phi_0, m_0)\|_{H^1}^2+\|\Delta \phi_0\|_{L^2}^2\}\nonumber\\
&&\quad +C\kappa_0\|m_0+\nabla \phi_0\|_{{\cal H}^1}^2. \label{Delta-w-ketusron}
\end{eqnarray}
\par
\parindent0pt
$\quad$
For the estimation of the $L^2-$norm on space of $w$
and the one on time-space of its spatial gradient,
let
\begin{eqnarray}
v = \int_0^t w (s) d s
\end{eqnarray}
and proceed to the energy estimates of $v$. 
From the direct calculation,
$v$ is seen to satisfy
\begin{eqnarray}
\left\{
\begin{array}{ll}
v_{t t} - \Delta v - \Delta v_t +\kappa_0 \Delta^2 v
=
\phi_0 + \Delta \mbox{div} \Phi,
\label{2-integral-eqn}
\\
v (0) = 0,
\qquad
v_t (0) = -\mbox{div} \Phi.
\end{array}
\right.
\end{eqnarray}
A test function
$v_t$ in (\ref{2-integral-eqn})
gives
\begin{eqnarray}
&&\frac{1}{2}
\frac{d}{d t}
\|v_t (t)\|_{L^2}^2
+
\frac{1}{2}
\frac{d}{d t}
\|\nabla v (t)\|_{L^2}^2
+
\|\nabla v_t (t)\|_{L^2}^2 +\frac{1}{2}\kappa_0
\frac{d}{d t}
\|\Delta v (t)\|_{L^2}^2\nonumber
\\
&&\qquad =
(\phi_0 + \Delta \mbox{div} \Phi, v_t (t))
\nonumber
\\
&&\qquad =
\frac{d}{d t}
(\phi_0 + \Delta \mbox{div} \Phi, v (t)),
\label{0v-l2-est}
\end{eqnarray}
being integrated on $(0, t)$
and yielding
\begin{eqnarray}
&&
\frac{1}{2}
\|v_t (t)\|_{L^2}^2
+
\frac{1}{2}
\|\nabla v (t)\|^2
+
\int_0^t
\|\nabla v_s (s)\|_{L^2}^2
d s
+
\frac{1}{2}
\kappa_0\|\Delta v (s)\|_{L^2}^2 
\nonumber
\\
&&
=
(\phi_0 + \Delta \mbox{div} \Phi, v (t))
+\frac{1}{2}\|w(0)\|_{L^2}^2.
\label{1v-L2-est}
\end{eqnarray}
Now the 1st term of the right hand side of (\ref{1v-L2-est}) is
%
\begin{eqnarray}
\left(
\phi_0 + \Delta \mbox{div} \Phi, v (t)
\right)
=
\left(
\phi_0, v (t)
\right)
+
\left(
\Delta \mbox{div} \Phi, v (t)
\right)
\label{2v-L2-est}
\end{eqnarray}
of which the 1st term is evaluated
by the use of Lemma \ref{Fefferman-Stein inequality}, \ref{poancare}
and Young's inequality as
\begin{eqnarray}
\left|
\left(
\phi_0, v (t)
\right)
\right|
\leq
C \|\phi_0\|_{{\cal H}^1}
\|v (t)\|_{{\rm BMO}}
&\leq&
C \|\phi_0\|_{{\cal H}^1}
\|\nabla v (t)\|_{L^2}
\nonumber
\\
&\leq&
\frac{1}{8}
\|\nabla v (t)\|_{L^2}^2
+
C \|\phi_0\|_{{\cal H}^1}^2
\nonumber
\end{eqnarray}
and, the 2nd term is controlled by \eqref{Phi-est1} as
\begin{eqnarray}
\left|\left(
 \Delta \mbox{div} \Phi, v (t)
\right)
\right|
&=&
\left|\left(
\Delta \Phi,  \nabla v (t)
\right)
\right|
\nonumber
\\
&\leq&
\|m_0+ \nabla \phi_0\|_{{\cal H}^1}
\|\nabla v (t)\|_{L^2}
\nonumber
\\
&\leq&
\frac{1}{8}
\|\nabla v (t)\|_{L^2}^2
+
C\|m_0+ \nabla \phi_0\|_{{\cal H}^1}^2.
\nonumber
\end{eqnarray}
Since $w (0) = - \mbox{div} \Phi$,
the right hand side of (\ref{1v-L2-est}) is bounded by
\begin{eqnarray}
&&(\phi_0 + \Delta \mbox{div} \Phi, v (t))
+\frac{1}{2}\|w(0)\|_{L^2}^2\nonumber\\
&&
\quad
\leq
\frac{1}{4}
\|\nabla v (t)\|_{L^2}^2
+
C 
\|\phi_0\|_{{\cal H}^1}^2
+
C 
\|m_0+ \nabla \phi_0\|_{{\cal H}^1}^2.
\nonumber
\end{eqnarray}
and thus, it follows that
\begin{eqnarray}
&&
\frac{1}{2}
\|v_t (t)\|^2
+
\frac{1}{4}
\|\nabla v (t)\|^2
+
\int_0^t
\|\nabla v_s (s)\|^2
\,
d s
+
\frac{1}{2}
\kappa_0\|\Delta v (s)\|_{L^2}^2 
\nonumber
\\
&&
\leq
C 
\|\phi_0\|_{{\cal H}^1}^2
+
C 
\|m_0+ \nabla \phi_0\|_{{\cal H}^1}^2.
\label{3v-L2-est}
\end{eqnarray} 
Gathering \eqref{w-base}, \eqref{nabla-w-est}, \eqref{Delta-w-ketusron} and \eqref{3v-L2-est}, we obtain Theorem \ref{time-spatial est}. 
\par
\vskip5pt
\parindent0pt
$\quad$
Based on Theorem \ref{time-spatial est}, Theorem \ref{time-decay est-CNSK} is proved as follows. We set a total energy of $w$ as 
$$
E(w(t))=\|w_t(t)\|_{L^2}^2+\|\nabla w(t)\|_{L^2}^2+\kappa_0 \|\Delta w(t)\|_{L^2}^2.
$$
By the proof of Theorem \ref{time-spatial est}
and integration by parts, we find that 
\begin{eqnarray}
	CJ_0
	&\geq&
	\int^t_0
	E (w (s))
	d s
	\nonumber
	\\
	&=&
	t E (w (t))
	- \int^t_0
	s
	\frac{d}{d s} E (w (s))
	d s.
	\nonumber
\end{eqnarray}
Since, by integration by parts again,
\begin{eqnarray}
	\frac{d}{d t} E (w (t))
	&=&
	\frac{d}{dt}\Big\{
	\|\nabla w (t)\|_{L^2}^2
	+
	\|w_t (t)\|_{L^2}^2
	+
	\kappa_0\|\Delta w(t)\|_{L^2}^2
	\Big\}
	\nonumber
	\\
	&=&
	2(\nabla w_t (t), \nabla w (t))
	+
	2(w_{tt} (t), w_t (t))
	+
	2(\kappa_0 \Delta w, \Delta w_t)
	\nonumber
	\\
	&=&
	-2 (w_t (t), \Delta w (t))
	+ 2(w_{tt} (t), w_t (t))
	+
	2(\kappa_0 \Delta^2 w,  w_t)
	\nonumber
	\\
	&=&
	2(w_{tt} (t) - \Delta w (t) +\kappa_0 \Delta^2 w(t), 
	w_t (t))
	\nonumber
	\\
	&=&
	2(\Delta w_t (t), w_t (t))
	\nonumber
	\\
	&=&
	- 2\|\nabla w_t (t)\|_{L^2}^2
	\nonumber
	\\
	&=&
	-2 \|\nabla \phi (t)\|_{L^2}^2,
	\nonumber
\end{eqnarray}
we have
\begin{eqnarray}
	CJ_0
	&\geq&
	t E (w (t))
	+
	2\int^t_0
	s
	\|\nabla \phi(s)\|^2
	d s.
	\label{B estim}
\end{eqnarray}
This together with \eqref{nabla-w-est} gives the assertion of Theorem \ref{time-decay est-CNSK}. The proof is completed. $\hfill\square$

\subsection{Proof of Theorem \ref{time-spatial est-diffusive2}}

In this section, we show the validity of Theorem \ref{time-spatial est-diffusive2}.  
By taking the Fourier transform of \eqref{cnsk-linear} with respect to the space variable $x$,
we have the following ordinary differential equation with a parameter $\xi$. 
\begin{eqnarray}
\left\{
\begin{array}{lll}
\partial_{t}\hat{\phi}(t,\xi) +i\gamma\xi\cdot\hat{m}(t,\xi)=0,\\
\partial_{t}\hat{m}(t,\xi)+\nu|\xi|^2 \hat{m}(t,\xi)+\tilde{\nu}\xi(\xi\cdot \hat{m}(t,\xi))+i\gamma \xi \hat{\phi}(t,\xi)
+i\xi \kappa_0 |\xi|^2\hat{\phi}(t,\xi)=0,\\
\hat{\phi}(0,\xi)=\hat{\phi}_0, \ \ \hat{m}(0,\xi)=\hat{m}_0.
\label{cnsk-f}
\end{array}
\right.
\end{eqnarray}
Therefore, the solutions of \eqref{cnsk-linear} are given by the following formulas by \cite{Kobayashi-Tsuda}. 
Let $K=\displaystyle\frac{2\sqrt{\kappa_0\gamma}}{\nu+\tilde{\nu}}$ and $B=\displaystyle\frac{2\gamma}{\nu+\tilde{\nu}}$. 
For $|\xi|\neq 0,  B/\sqrt{1-K^2}$ when $0< K < 1$ and $|\xi|\neq 0$ when $K \geq 1$,
the Fourier transforms of $\phi$ and $m$ are given explicitly by the formulas; 
\begin{eqnarray}
\hat{\phi}&=&\displaystyle\frac{\lambda_+(\xi) e^{\lambda_-(\xi) t}-\lambda_-(\xi) e^{\lambda_+(\xi) t}}
{\lambda_+(\xi) -\lambda_-(\xi)}\hat{\phi}_0-i\gamma \displaystyle\frac{ e^{\lambda_+(\xi) t}-e^{\lambda_-(\xi) t}}
{\lambda_+(\xi) -\lambda_-(\xi)}\xi\cdot \hat{m}_0,\nonumber\\
\hat{m}&=&e^{-\nu|\xi|^2 t}\hat{m}_0-i\xi(\gamma+\kappa_0|\xi|^2)\left(\displaystyle\frac{ e^{\lambda_+(\xi) t}-e^{\lambda_-(\xi) t}}
{\lambda_+(\xi) -\lambda_-(\xi)}\right)\hat{\phi}_0\nonumber\\
&&\quad +\left(\displaystyle\frac{\lambda_+(\xi) e^{\lambda_+(\xi) t}-\lambda_-(\xi) e^{\lambda_-(\xi) t}}{\lambda_+(\xi) -\lambda_-(\xi)}-e^{-\nu|\xi|^2t}\right)\frac{\xi(\xi\cdot\hat{m}_0)}{|\xi|^2},\label{solution-linear}
\end{eqnarray}
where $\lambda_\pm$ are given by 
\begin{eqnarray} 
\lambda_{\pm}(\xi)=-A(|\xi|^2\pm\sqrt{|\xi|^4-B^2|\xi|^2-K^2|\xi|^4})\label{chracter}
\end{eqnarray}
with a positive constant $A=\displaystyle\frac{\nu+\tilde{\nu}}{2}$  
and stand for roots of the characteristic equation of \eqref{cnsk-f}. 
If $0<K< 1$ and $\min\Big\{\frac{1}{2}, \, \displaystyle\frac{B}{2\sqrt{1-K^2}}\Big\}\leq |\xi| \leq 2\displaystyle\frac{B}{\sqrt{1-K^2}}$ 
$\hat{\phi}$ and  $\hat{m}$ are represented as 
\begin{eqnarray}
\hat{\phi}&=&\displaystyle\frac{1}{2\pi i}\oint_{\Gamma}\frac{(z+A|\xi|^2)e^{zt}}
{z^2+(\nu+\tilde{\nu})|\xi|^2 z+\kappa_0\gamma|\xi|^4+\gamma^2|\xi|^2}dz\hat{\phi}_0\nonumber\\
&&\quad-\displaystyle\frac{\gamma}{2\pi }\oint_{\Gamma}\frac{e^{zt}}
{z^2+(\nu+\tilde{\nu})|\xi|^2 z+\kappa_0\gamma|\xi|^4+\gamma^2|\xi|^2}dz\xi\cdot\hat{m}_0,\label{solmidlerho}\\
\hat{m}&=&e^{-\nu|\xi|^2 t}\hat{m}_0-\displaystyle\frac{\gamma \xi}{2\pi }
\oint_{\Gamma}\frac{e^{zt}}
{z^2+(\nu+\tilde{\nu})|\xi|^2 z+\kappa_0\gamma|\xi|^4+\gamma^2|\xi|^2}dz\hat{\phi}_0\nonumber\\
&&\quad+\left(\displaystyle\frac{1}{2\pi i}\oint_{\Gamma}\frac{ze^{zt}}
{z^2+(\nu+\tilde{\nu})|\xi|^2 z+\kappa_0\gamma|\xi|^4+\gamma^2|\xi|^2}dz\right)\frac{\xi(\xi\cdot\hat{m}_0)}{|\xi|^2},\label{solmidlem}
\end{eqnarray}
where $\Gamma$ is a closed pass surrounding
$\lambda_{\pm}$ and included in the set $\{z\in \mathbb{C}|{\rm Re}z\leq -c_0\}$
and $c_0$ is
a positive number satisfying
$$
\max_{\min\Big\{\frac{1}{2}, \,\frac{B}{2\sqrt{1-K^2}}\Big\}\leq|\xi|\leq 2\frac{B}{\sqrt{1-K^2}}}{\rm Re}\lambda_{\pm}\leq -2c_0.
$$
\par
\parindent0pt
$\quad$
Cut-off functions $\varphi_1$,
$\varphi_\infty$ and $\varphi_M$ in $C^{\infty}(\mathbb{R}^2)$ are defined by \cite{Kobayashi-Tsuda} as follows:
in the case such that $K \neq 1$,
$\varphi_1$ is given by 
\begin{eqnarray*}
&&
\varphi_1(\xi)=
\left\{
\begin{array}{ll}
1 \ \ \mbox{for} \ \   |\xi|\leq \frac{1}{2}
\\
0 \ \ \mbox{for} \ \  |\xi|\geq 1
\end{array}
\right.
\quad
\mbox{if} \ \frac{B}{2\sqrt{|1-K^2|}}> 1 ; 
\\ 
&&
\varphi_1(\xi)=
\left\{
\begin{array}{ll}
1 \ \ \mbox{for} \ \   |\xi|\leq \displaystyle\frac{B}{2\sqrt{|1-K^2|}}
\\
0 \ \ \mbox{for} \ \  |\xi|\geq 1
\end{array}
\right.
\quad
\mbox{if} \ \frac{B}{2\sqrt{|1-K^2|}}\leq 1 < \frac{B}{\sqrt{2|1-K^2|}} ;
\\
&&
\varphi_1(\xi)=
\left\{
\begin{array}{ll}
1 \ \ \mbox{for} \ \   |\xi|\leq \displaystyle\frac{B}{2\sqrt{|1-K^2|}}
\\
0 \ \ \mbox{for} \ \  |\xi|\geq \displaystyle\frac{B}{\sqrt{2|1-K^2|}}
\end{array}
\right.
\quad
\mbox{if} \ \frac{B}{\sqrt{2|1-K^2|}}\leq 1.
\end{eqnarray*}
$\varphi_\infty$ and $\varphi_M$ are
\begin{eqnarray*}
&&
\varphi_\infty(\xi)=
\left\{
\begin{array}{ll}
1 \ \  \mbox{for} \ \ |\xi|\geq \displaystyle\frac{2B}{\sqrt{|1-K^2|}}
\\
0 \ \  \mbox{for} \ \ |\xi|\leq \displaystyle\frac{\sqrt{2}B}{\sqrt{|1-K^2|}},
\end{array}
\right.
\quad
\varphi_M(\xi)=1-\varphi_1(\xi)-\varphi_\infty(\xi).
\end{eqnarray*}
In the case that $K=1$,
$\varphi_1$ and $\varphi_\infty$ are
\begin{eqnarray*}
&&
\varphi_1(\xi)=
\left\{
\begin{array}{ll}
1 \ \ \mbox{for} \ \   |\xi|\leq \frac{1}{2}
\\
0 \ \ \mbox{for} \ \  |\xi|\geq 1,
\end{array}
\right. 
\quad
\varphi_\infty(\xi)=
\left\{
\begin{array}{ll}
0 \ \ \mbox{for} \ \   |\xi|\leq \frac{1}{2}
\\
1 \ \ \mbox{for} \ \  |\xi|\geq 1,
\end{array}
\right.
\\
&&
\varphi_1(\xi)+\varphi_\infty(\xi)=1. 
\end{eqnarray*}
We define the solution operators $E_1$ and $E_\infty$ on a low frequency part and on a high frequency part of \eqref{cnsk-linear}, respectively, as follows :
\begin{eqnarray}
&&E_1(t)=(E_{1,\phi}(t),E_{1,m}(t)),
\nonumber
\\
&&E_\infty(t)=(E_{\infty,\phi}(t),E_{\infty,m}(t)),
\label{sol-ope-high}
\end{eqnarray}
where
\begin{eqnarray}
&&E_{1,\phi}(t)(\phi_0,m_0)(x)=\mathcal{F}^{-1}[\varphi_1(\xi)\hat{\phi}(t,\xi)](x),\label{sol-ope-low1}\\
&&E_{1,m}(t)(\phi_0,m_0)(x)=\mathcal{F}^{-1}[\varphi_1(\xi)\hat{m}(t,\xi)](x),\label{sol-ope-low2}\\
&&E_{\infty,\phi}(t)(\phi_0,m_0)(x)=\mathcal{F}^{-1}[(\varphi_M(\xi)+\varphi_\infty(\xi))\hat{\phi}(t,\xi)](x),\label{sol-ope-high1}\\
&&E_{\infty,m}(t)(\phi_0,m_0)(x)=\mathcal{F}^{-1}[(\varphi_M(\xi)+\varphi_\infty(\xi))\hat{m}(t,\xi)](x).\label{sol-ope-high2}
\end{eqnarray}
In \cite{Kobayashi-Tsuda} the solution operator is shown to
have an exponential decay in time on the high frequency part \eqref{sol-ope-high}.
In fact we have

\vspace{2ex}

\begin{thm}\label{est-linear-high}{\rm \cite[Theorem 3.2]{Kobayashi-Tsuda}} 
Let $1\leq p \leq \infty$. Then it holds that
\begin{eqnarray}
\lefteqn{\|\del_t^k \del_x^\alpha E_{\infty}(t)(\phi_0,m_0)\|_{L^p}}\nonumber\\
&\leq & C_{k,\alpha}e^{-ct} \Big\{(1+t^{-\delta_1-\frac{|\alpha|}{2}-k})\Big[\|\phi_0\|_{L^p}+\|m_0\|_{L^p}\Big]+
t^{-\delta_2 -\frac{|\alpha|}{2}-k}\|\phi_0\|_{L^p}\Big\}\label{est-linear-high1}
\end{eqnarray}
for $t>0$, $k \geq 0$ and $|\alpha|\geq 0$,  
where $(\delta_1,\delta_2 )=(1/2,1)$ and $(1,3/2)$ for $K\neq 1$ and $K=1$, respectively.
\end{thm}

\vspace{2ex}

Therefore, in order to show Theorem \ref{time-spatial est-diffusive2}, it is enough to consider the low frequency part. We estimate the Green function.  
We put 
\begin{eqnarray*}
L_{11,j}(t,x)&=&\mathcal{F}^{- 1}\Big\{\displaystyle\frac{\lambda_+(\xi) e^{\lambda_-(\xi) t}-\lambda_-(\xi) e^{\lambda_+(\xi) t}}{\lambda_+(\xi) -\lambda_-(\xi)}\varphi_j(\xi)\Big\}(x), \\
L_{12,j}(t,x)&=&\mathcal{F}^{- 1}(-i\gamma\trans\xi\hat{L}_j),\\
\hat{L}_j(t,\xi)&=&\displaystyle\frac{ e^{\lambda_+(\xi) t}-e^{\lambda_-(\xi) t}}{\lambda_+(\xi) -\lambda_-(\xi)}\varphi_j(\xi), \\
L_{21,j}(t,x)&=&\mathcal{F}^{-1}\{-\xi(i\gamma +k_0|\xi|^2 )\hat{L}_j\},\\
L_{22,j}(t,x)&=&K_{1,j}(t,x)+K_{2,j}(t,x)-K_{3,j}(t,x),\\
K_{1,j}(t,x)&=&\mathcal{F}^{- 1}\Big[e^{-\nu|\xi|^2 t}\varphi_j(\xi)\Big](x)I_n,\\
K_{2,j}(t,x)&=&\mathcal{F}^{- 1}\Big\{\displaystyle\frac{\lambda_+(\xi) e^{\lambda_+(\xi) t}-\lambda_-(\xi) e^{\lambda_-(\xi) t}}{\lambda_+(\xi) -\lambda_-(\xi)}\frac{\xi \trans{\xi}}{|\xi|^2}\varphi_j(\xi)\Big\}(x), \\
K_{3,j}(t,x)&=&\mathcal{F}^{- 1}\Big[e^{-\nu|\xi|^2 t}\frac{\xi \trans{\xi}}{|\xi|^2}\varphi_j(\xi)\Big](x)
\end{eqnarray*}
for $j=1,\infty$. We see from $\eqref{solution-linear}$ that 
\begin{eqnarray}
E_j(t)(\phi_0,m_0)=
\begin{pmatrix}
L_{11,j}(t,\cdot) & L_{12,j}(t,\cdot)\\
L_{21,j}(t,\cdot) & L_{22,j}(t,\cdot)\\
\end{pmatrix}
*
\begin{pmatrix}
\phi_0\\
m_0
\end{pmatrix}
\label{convolution-solutionfomula}
\end{eqnarray}
for $j=1, \infty$. 

We set  
$$
K_{1}m_0=\mathcal{F}^{- 1}\Big[\frac{\lambda_{+}(\xi)e^{\lambda_+(\xi) t}-\lambda_-(\xi)e^{\lambda_-(\xi) t}}{\lambda_+(\xi) -\lambda_-(\xi)}\frac{\xi \trans{\xi}}{|\xi|^2}\varphi_1(\xi)\hat{m}_0\Big]. 
$$
Note that $K_1 m_0$ is a part of the Green matrix and corresponds to the diffusion wave part $m-{\cal K}_{\nu}\ast m_{0,in}$. Our claim is the following estimate

\vspace{2ex}

\begin{prop}\label{low-kernelest}
Let $m_0 \in L^2 \cap L^1$, $|x| m_0 \in L^1$ and $\hat{m}_0 \neq 0$. Then it holds that 
$$
\limsup_{t\rightarrow \infty}\Big|\displaystyle\int_0^t \|K_1 m_0 (\tau)\|_{L^2}^2 d\tau -\log t \Big| \leq C_1,
$$
\end{prop}
where $C_1=C_1(u_0)$ is a positive constant independent of $t$. 

\vspace{2ex}

{\bf Proof.}  By the Plancherel theorem and \eqref{chracter}, we see that there exists a positive constant $C$ such that 
$$
C\|e^{-|\xi|^2 t} \hat{m_0}(\xi)\|_{L^2 (|\xi| \leq c_1)}^2\leq \|K_1 m_0 (t)\|_{L^2}^2 \leq \frac{1}{C}\|e^{-|\xi|^2 t} \hat{m_0}(\xi)\|_{L^2 (|\xi| \leq c_1)}^2, 
$$
where $c_1=1$ when $\frac{B}{2\sqrt{|1-K^2|}} >1$ or $\frac{B}{2\sqrt{|1-K^2|}} \leq 1 < \frac{B}{\sqrt{2|1-K^2|}}$ and $c_1=\frac{B}{\sqrt{2|1-K^2|}}$ when $\frac{B}{\sqrt{2|1-K^2|}} \leq 1$.
%
%
Hence we have to estimate $\|e^{-|\xi|^2 t} \hat{m_0}(\xi)\|_{L^2 (|\xi| \leq c_1)}^2$.
It follows from the polar coordinates that 
\begin{eqnarray*}
\|e^{-|\xi|^2 t} \hat{m_0}(\xi)\|_{L^2 (|\xi| \leq c_1)}^2=C_\omega \displaystyle\int_{0}^{c_1} |e^{-r^2 t} \hat{m}_0(r\omega)|^2 r dr,
\end{eqnarray*}
where $r=|\xi|$, $\omega=\xi/|\xi|$ and $C_{\omega}$ 
is some positive constant which appears in the polar coordinates.
Changing variables $r\sqrt{t}=s$, we have 
$$
\displaystyle\int_{0}^{c_1} |e^{-r^2 t} \hat{m}_0(\omega r)|^2 r dr=t^{-1}\displaystyle\int_{0}^{c_1\sqrt{t}} e^{-2s^2} |\hat{m_0}(s\omega/\sqrt{t})|^2 s ds.
$$
This together with  the fundamental theorem of calculus for $\hat{m}_0$ imply that 
\begin{eqnarray}
\lefteqn{\Big|\|e^{-|\xi|^2 t} \hat{m_0}(\xi)\|_{L^2 (|\xi| \leq c_1)}^2-C_\omega t^{-1}|\hat{m}_0(0)|^2\displaystyle
\int_{0}^{c_1\sqrt{t}} e^{-2s^2}  s ds\Big|}
\nonumber
\\
&& =C_{\omega}\Big|t^{-1}\displaystyle\int_{0}^{c_1\sqrt{t}}
(|\hat{m}_0(s\omega/\sqrt{t})|^2-|\hat{m}_0(0)|^2) e^{-2s^2}  s ds\Big|
\nonumber
\\
&& \leq Ct^{-1}\displaystyle\int_{0}^{c_1\sqrt{t}}
|\hat{m}_0(s\omega/\sqrt{t})-\hat{m}_0(0)|^2e^{-2s^2}  s ds
%
%
\nonumber
\\
&&=Ct^{-1}\displaystyle\int_{0}^{c_1\sqrt{t}}
\Big|\displaystyle\int_{0}^{1}\nabla_{\xi}\hat{m}_0\Big(\theta \frac{s\omega}{\sqrt{t}}\Big)d\theta \frac{s\omega}{\sqrt{t}}\Big|^2e^{-2s^2}  s ds
\nonumber
\\
&& \leq Ct^{-1}\displaystyle\int_{0}^{c_1\sqrt{t}}
\|\nabla_{\xi}\hat{m}_0\|_{L^\infty}^2\Big|\frac{s\omega}{\sqrt{t}}\Big |^2 e^{-2s^2}  s ds
\nonumber
\\
&& \leq C\||x| m_0\|_{L^1}^2 t^{-2}\displaystyle\int_{0}^{c_1\sqrt{t}}e^{-2s^2}  s^3 ds \nonumber\\
&& \leq C_2 t^{-2}\label{est-2-1}
\end{eqnarray}
for a positive constant $C_2=C_2(u_0)$ independent of $t$. Note that $\hat{m}_0 (0) \neq 0$ by our assumption.  Since 
 $$
\displaystyle\int_{0}^{1}\|e^{-|\xi|^2 \tau} \hat{m_0}(\xi)\|_{L^2}^2 d\tau \leq C\|m_0\|_{L^2}^2, 
 $$
we estimate $\displaystyle\int_{1}^{t}\|e^{-|\xi|^2 \tau} \hat{m_0}(\xi)\|_{L^2}^2 d\tau$. We set 
$$
I_1(t):=C_\omega|\hat{m}_0(0)|^2 t^{-1}\displaystyle\int_{0}^{c_1\sqrt{t}}e^{-2s^2}  s ds
$$
Applying \eqref{est-2-1} yields that for $t\geq 1$ 
\begin{eqnarray}
\Big|\displaystyle\int_{1}^{t} \Big\{ \|e^{-|\xi|^2 \tau} \hat{m_0}(\xi)\|_{L^2 (|\xi| \leq c_1)}^2 - I_1(\tau)\Big\}
d\tau\Big| \leq C_2\displaystyle\int_{1}^{t} \tau^{-2}d\tau \leq C_2. \label{keturon1}
\end{eqnarray}
Let positive constants $Q_1$ and $Q_2$ be defined by
$$
Q_1:=\displaystyle\int_{0}^{c_1}e^{-2s^2}  s^3 ds,  \ \ Q_2:=\displaystyle\int_{0}^{\infty}e^{-2s^2}  s^3 ds.
$$
It obviously follows that for $t \geq 1$
$$
Q_1 \leq \displaystyle\int_{0}^{c_1\sqrt{t}}e^{-2s^2}  s^3 ds \leq Q_2.
$$
This implies that for $t \geq 1$
$$
 \displaystyle\int_{1}^{t}C_\omega|\hat{m}_0(0)|^2 Q_1 \tau^{-1}d\tau \leq \displaystyle\int_{1}^{t} I_1(\tau)d\tau \leq  \displaystyle\int_{1}^{t}C_\omega|\hat{m}_0(0)|^2 Q_2 \tau^{-1}d\tau
$$
and thus, 
\begin{eqnarray}
C_\omega|\hat{m}_0(0)|^2 Q_1 \log t \leq \displaystyle\int_{1}^{t} I_1(\tau)d\tau \leq C_\omega|\hat{m}_0(0)|^2 Q_2\log t. \label{keturon2}
\end{eqnarray}
\eqref{keturon1} and \eqref{keturon2} yield that 
\begin{eqnarray*}
\limsup_{t \rightarrow \infty}\Big|\displaystyle\int_{1}^{t} \|e^{-|\xi|^2 \tau} \hat{m_0}(\xi)\|_{L^2 (|\xi| \leq c_1)}^2 d\tau - \log t\Big|
\leq 
CC_2. 
\end{eqnarray*}
Therefore, there exists a positive constant $C_1=C_1(u_0)$ 
independent of $t$ such that
$$
\limsup_{t\rightarrow \infty}\Big|\displaystyle\int_0^t \|K_1 m_0 (\tau)\|_{L^2}^2 d\tau -\log t \Big| \leq C_1. 
$$

\vspace{2ex}
\par
\parindent0pt
$\quad$
Since other parts of the diffusion wave $m-{\cal K}_{\nu}\ast m_{0,in}$ which appear in the Green matrix on the low frequency part are estimated similarly to  Proposition \ref{low-kernelest}, we obtain the estimation in Theorem \ref{time-spatial est-diffusive2}.
The proof is completed. $\hfill\square$

%
%
\section{Appendix}
Here we will demonstrate the proof of Theorem \ref{exist-d-eqn}. 
\par
\parindent0pt
{\it Proof of Theorem \ref{exist-d-eqn}}.
Now we define an operator $T$ for $f \in {\cal H}^1$ by 
\begin{eqnarray}
T f = {\cal F}^{-1}((1+\kappa_0 |\xi|^2)^{-1} \hat{f}):=K\ast f.
\nonumber
\end{eqnarray}
From direct computation we see that 
\begin{eqnarray}
|\del_{\xi}^{\alpha} \hat{K}(\xi)| \leq C|\xi|^{-|\alpha|},
\quad
\mbox{for any $\xi \neq 0$
and $|\alpha|\geq 0$}
\nonumber
\end{eqnarray}
for a positive constant $C$ independent of $\kappa_0$.
Then, it follows from Shimizu and Shibata \cite[Theorem 2.3]{Shibata-Shimizu} that 
$$
|\del_{\xi}^{\alpha}K(x)| \leq C_1|x|^{-2-|\alpha|} \ \ (x \neq 0)
$$
holds true for a positive constant $C_1$ independent of $\kappa_0$. 
By this fact and the multiplier type theorem on the Hardy space as in Stein \cite[Chapter 3, Section 3.2, Theorem 4]{Stein}, we find that $T$ is a bounded operator on ${\cal H}^1$ and 
 \begin{eqnarray}
 \|T f \|_{{\cal H}^1} \leq C \|f\|_{{\cal H}^1}, 
 \nonumber
 \end{eqnarray}
 where $C$ is independent of $\kappa_0$.
Therefore, $(I - \kappa_0 \Delta)^{- 1}$ is bounded on ${\cal H}^1$
and thus, $(I - \kappa_0 \Delta)^{- 1} (m_0 + \nabla \phi_0) \in {\cal H}^1$
for $m_0 + \nabla \phi_0 \in {\cal H}^1$.
Then, from \cite{Coifman} we obtain the existence of a solution to (\ref{doubly-eqn})
satisfying (\ref{Phi-est1}).
Furthermore, since $(1-\kappa_0 \Delta)^{-1}$ is bounded on $L^2$ we have the estimation
\begin{eqnarray}
\|\div \Delta \Phi\|_{L^2}
&=&
\|\div (1-\kappa_0 \Delta)^{-1}(m_0 + \nabla \phi_0)\|_{L^2}
\nonumber
\\
&\leq&
C\|\div (m_0 + \nabla \phi_0)\|_{L^2}
\nonumber
\\
&\leq&
C(\|\Delta \phi_0\|_{L^2} + \|\trans(\phi_0, m_0)\|_{H^1}),
\nonumber
\end{eqnarray}
from which (\ref{Phi-est2}) is obtained.
The proof is completed.
\qed

\vspace{2ex}

\noindent {\bf Acknowledgments.} 
The first author is partly supported by Grants-in-Aid for Scientific Research with the Grant number: 16H03945.
The second author is partly supported by
Grants-in-Aid for Scientific Research
with the number: 18K03375.
The third author is partly supported by Grant-in-Aid for JSPS Fellows with the Grant number: A17J047780.


\end{document}